\documentclass{llncs}[10pt]
\usepackage{llncsdoc}
\usepackage{amsmath,graphicx}
\usepackage{amsfonts} 
\usepackage{amssymb}
\usepackage[latin1]{inputenc}

\title{Riemannian Gaussian distributions on the space of positive-definite quaternion matrices}

\author{Salem Said$^{\,1}$, Nicolas Le Bihan$^{\,2}$, Jonathan H. Manton$^{\,3}$}

\institute{1. Laboratoire IMS (CNRS - UMR 5218), 2. Gipsa-lab (CNRS - UMR 5216), \\ 3. The University of Melbourne, Dept. of Electrical and Electronic Engineering}

\usepackage[update,prepend]{epstopdf}
\RequirePackage{fix-cm}
\DeclareMathSizes{10}{10}{5}{5}

\begin{document}

\maketitle

\begin{abstract}
 Recently, Riemannian Gaussian distributions were defined on spaces of positive-definite real and complex matrices. The present paper extends this definition to the space of positive-definite quaternion matrices. In order to do so, it develops the Riemannian geometry of the space of positive-definite quaternion matrices, which is shown to be a Riemannian symmetric space of non-positive curvature. The paper gives original formulae for the Riemannian metric of this space, its geodesics, and distance function. Then, it develops the theory of Riemannian Gaussian distributions, including the exact expression of their probability density, their sampling algorithm and statistical inference. 
\end{abstract}
\begin{keywords}
Riemannian Gaussian distribution, quaternion, positive-definite \\matrix, symplectic group, Riemannian barycentre
\end{keywords}
\section{Introduction} The Riemannian geometry of the spaces $\mathcal{P}_n$ and $\mathcal{H}_n\,$, respectively of $n \times n$ positive-definite real and complex matrices, is well-known to the information science community~\cite{pennec,moakher}. These spaces have the property of being Riemannian symmetric spaces of non-positive curvature~\cite{helgasson,besse},
$$
\mathcal{P}_n = \left.\mathrm{GL}(n,\mathbb{R})\middle/\mathrm{O}(n)\right. \hspace{0.5cm}
\mathcal{H}_n = \left.\mathrm{GL}(n,\mathbb{C})\middle/\mathrm{U}(n)\right.
$$
where $\mathrm{GL}(n,\mathbb{R})$ and $\mathrm{GL}(n,\mathbb{C})$ denote the real and complex linear groups, and $\mathrm{O}(n)$ and $\mathrm{U}(n)$ the orthogonal and unitary groups. Using this property, Riemannian Gaussian distributions were recently introduced on $\mathcal{P}_n$ and $\mathcal{H}_n$~\cite{said1,said2}. The present paper introduces the Riemannian geometry of the space $\mathcal{Q}_n$ of $n \times n$ positive-definite quaternion matrices, which is also a Riemannian symmetric space of non-positive curvature~\cite{besse},
$$
\mathcal{Q}_n \,=\, \left. \mathrm{GL}(n,\mathbb{H})\middle/\mathrm{Sp}(n) \right.
$$
where $\mathrm{GL}(n,\mathbb{H})$ denotes the quaternion linear group, and $\mathrm{Sp}(n)$ the compact symplectic group. It then studies Riemannian Gaussian distributions on $\mathcal{Q}_n$. The main results are the following\,: Proposition \ref{prop:metric} gives the Riemannian metric of the space $\mathcal{Q}_n$, Proposition \ref{prop:polar1} expresses this metric in terms of polar coordinates on the space $\mathcal{Q}_n$, Proposition \ref{prop:polar2} uses Proposition \ref{prop:polar1} to compute the moment generating function of a Riemannian Gaussian distribution on $\mathcal{Q}_n$, and Propositions \ref{prop:sampling} and \ref{prop:inference} describe the sampling algorithm and maximum likelihood estimation of Riemannian Gaussian distributions on $\mathcal{Q}_n$.  Motivation for studying matrices from $\mathcal{Q}_n$ comes from their potential use in multidimensional bivariate signal processing \cite{flamant2016}.


%
%
%
%

\section{Quaternion matrices, $\mathrm{GL}(\mathbb{H})$ and $\mathrm{Sp}(n)$}
Recall the non-commutative division algebra of quaternions, denoted $\mathbb{H}$, is made up of elements $q = q_0 \,+\, q_1\,\mathrm{i} +\, q_2\,\mathrm{j}+\, q_3\,\mathrm{k}$ where $q_0, q_1, q_2, q_3 \in \mathbb{R}$, and the imaginary units $\mathrm{i}, \mathrm{j}, \mathrm{k}$ satisfy the relations~\cite{conway}
\begin{equation} \label{eq:quaternions}
  \mathrm{i}^2 = \mathrm{j}^2 = \mathrm{k}^2 = \mathrm{ijk} \,=\, -1 
\end{equation}
The real part of $q$ is $\mathrm{Re}(q) = q_0\,$, its conjugate is $\bar{q} = q_0 \,-\, q_1\,\mathrm{i} -\, q_2\,\mathrm{j}-\, q_3\,\mathrm{k}$ and its squared norm is $|q|^2 = q\bar{q}\,$. The multiplicative inverse of $q \neq 0 $ is given by $q^{-1} = \bar{q}/|q|^2\,$. \\[0.1cm]
\indent The set $M_n(\mathbb{H})$ consists of $n\times n$ quaternion matrices $A$~\cite{zhang}. These are arrays $A = (A_{ij}\,;\, i,j = 1,\ldots, n)$ where $A_{ij} \in \mathbb{H}$. The product $C = AB$ of $A, B \in M_n(\mathbb{H})$ is the element of $M_n(\mathbb{H})$ with
\begin{equation} \label{eq:quatmatrix}
  C_{ij} \,=\, \sum^n_{l=1}\, A_{il}B_{lj}
\end{equation}
A quaternion matrix $A$ is said invertible if it has a multiplicative inverse $A^{-1}$ with $AA^{-1} = A^{-1}A = I$ where $I$ is the identity matrix.  The conjugate-transpose of $A$ is $A^{\dagger}$ which is a quaternion matrix with $A^{\dagger}_{ij} = \bar{A}_{ji\,}$. 

The rules for computing with quaternion matrices are quite different from the rules for computing with real or complex matrices~\cite{zhang}. For example, in general, $\mathrm{tr}(AB) \neq \mathrm{tr}(BA)$, and $(AB)^T \neq B^{\,T}A^T$ where $^T$ denotes the transpose. For the results in this paper, only the following rules are needed~\cite{zhang},
\begin{equation} \label{eq:quatrules}
   (AB)^{-1} \,=\, B^{-1}A^{-1} \hspace{0.3cm} (AB)^{\dagger} \,=\, B^{\dagger}A^{\dagger} \hspace{0.3cm} \mathrm{Re}\, \mathrm{tr}(AB) = \mathrm{Re}\, \mathrm{tr}(BA)
\end{equation}
$\mathrm{GL}(n,\mathbb{H})$ consists of the set of invertible quaternion matrices $A \in M_n(\mathbb{H})$. The subset of $A \in \mathrm{GL}(n,\mathbb{H})$ such that $A^{-1} = A^{\dagger}$ is denoted $\mathrm{Sp}(n) \subset \mathrm{GL}(n,\mathbb{H})$. 

It follows from (\ref{eq:quatrules}) that $\mathrm{GL}(n,\mathbb{H})$ and $\mathrm{Sp}(n)$ are groups under the operation of matrix multiplication, defined by (\ref{eq:quatmatrix}). However, one has more. Both these groups are real Lie groups. Usually, $\mathrm{GL}(n,\mathbb{H})$ is called the quaternion linear group, and $\mathrm{Sp}(n)$ the compact symplectic group. In fact, $\mathrm{Sp}(n)$ is a compact connected Lie subgroup of $\mathrm{GL}(n,\mathbb{H})$~\cite{kirillov}. 

The Lie algebras of these two Lie groups are given by
\begin{equation} \label{eq:liealgebras}
  \mathfrak{gl}(n,\mathbb{H}) = M_n(\mathbb{H}) \hspace{0.35cm} \mathfrak{sp}(n) = \left \lbrace X \in \mathfrak{gl}(n,\mathbb{H})\,\middle|\, X + X^{\dagger} = 0 \, \right \rbrace
\end{equation}
with the bracket operation $[X,Y] = XY - YX$. The Lie group exponential is identical to the quaternion matrix exponential
\begin{equation} \label{eq:exponential}
 \exp(X) = \sum_{m\geq 0} \frac{X^m}{m!} \hspace{1cm} X \in \mathfrak{gl}(n,\mathbb{H})
\end{equation}
For $A \in \mathrm{GL}(n,\mathbb{H})$ and $X \in \mathfrak{gl}(n,\mathbb{H})$, let $\mathrm{Ad}(A)\cdot X = AXA^{-1}\,$. Then, 
\begin{equation} \label{eq:adjoint}
 A\,\exp(X)\,A^{-1} \,=\, \exp\left( \,\mathrm{Ad}(A)\cdot X\,\right)
\end{equation}
as can be seen from (\ref{eq:exponential}).

%
%
%
%

\section{The space $\mathcal{Q}_n$ and its Riemannian metric}
The space $\mathcal{Q}_n$ consists of all quaternion matrices $S \in M_n(\mathbb{H})$ which verify $S = S^{\dagger}$ and
\begin{equation} \label{eq:qn}
\sum^n_{i,j=1}\, \bar{x}_i\,S_{ij}\,x_j \, > 0 \hspace{0.3cm} \text{for all non-zero } (x_1\,,\ldots,\, x_n) \, \in \mathbb{H}^n
\end{equation}
In other words, $\mathcal{Q}_n$ is the space of positive-definite quaternion matrices. Note that, due to the condition $S = S^{\dagger}$, the sum in (\ref{eq:qn}) is a real number. 

Define now the action of $\mathrm{GL}(n,\mathbb{H})$ on $\mathcal{Q}_n$ by $A\cdot S = ASA^{\dagger}$ for $A \in \mathrm{GL}(n,\mathbb{H})$ and $S \in \mathcal{Q}_n\,$. This is a left action, and is moreover transitive. Indeed~\cite{zhang}, each $S \in \mathcal{Q}_n$ can be diagonalized by some $K \in \mathrm{Sp}(n)$,
\begin{equation} \label{eq:spectral}
  S = K\,\exp(R)\, K^{-1} \,=\, \exp\left( \,\mathrm{Ad}(K)\cdot R\,\right) \hspace{0.2cm} ;\, R \,\,\text{ real diagonal matrix}
\end{equation}
where the second equality follows from (\ref{eq:adjoint}). Thus, each $S \in \mathcal{Q}_n$ can be written $S = AA^{\dagger}$ for some $A \in \mathrm{GL}(n,\mathbb{H})$,  which is the same as $S = A\cdot I$. 

For $A \in \mathrm{GL}(n,\mathbb{H})$, note that $A \cdot I = I$ iff $AA^{\dagger} = I$, which means that $A \in \mathrm{Sp}(n)$. Therefore, as a homogeneous space under the left action of $\mathrm{GL}(n,\mathbb{H})$,
\begin{equation} \label{eq:quotient}
  \mathcal{Q}_n \,=\, \left. \mathrm{GL}(n,\mathbb{H})\,\middle/\,\mathrm{Sp}(n) \right.
\end{equation}
The space $\mathcal{Q}_n$ is a real differentiable manifold. In fact, if $\mathfrak{p}_n$ is the real vector space of $X \in \mathfrak{gl}(n,\mathbb{H})$ such that $X = X^{\dagger}$, then it can be shown $\mathcal{Q}_n$ is an open subset of $\mathfrak{p}_n\,$. Therefore, $\mathcal{Q}_n$ is a manifold, and for each $S \in \mathcal{Q}_n$ the tangent space $T_S\mathcal{Q}_n$ may be identified with $\mathfrak{p}_n$. Moreover, $\mathcal{Q}_n$ can be equipped with a Riemannian metric as follows.

Define on $\mathfrak{gl}(n,\mathbb{H})$ the $\mathrm{Sp}(n)$-invariant scalar product
\begin{equation} \label{eq:scalproduct}
   \left \langle X\middle|Y\right \rangle  \,= \mathrm{Re}\,\mathrm{tr}(XY^{\dagger}) \hspace{0.5cm} X,Y\,\in \mathfrak{gl}(n,\mathbb{H})
\end{equation}
For $u, v$ in $T_S\mathcal{Q}_n \simeq \mathfrak{p}_n\,$, let 
\begin{equation} \label{eq:metric1}
   (u,v)_S \,= \left \langle\,(A^{-1}) u(A^{-1})^{\dagger}\,\middle|\,(A^{-1}) v(A^{-1})^{\dagger}\,\right \rangle
\end{equation}
where $A$ is any element of $\mathrm{GL}(n,\mathbb{H})$ such that $S = A\cdot I\,$. 
\begin{proposition}[Riemannian metric] \label{prop:metric} \\
(i)  For each $S \in \mathcal{Q}_n\,$,  formula (\ref{eq:metric1}) defines a scalar product on $T_S\mathcal{Q}_n \simeq \mathfrak{p}_n\,$, which is independent of the choice of $A$. \\
(ii) Moreover,  
\begin{equation} \label{eq:metric2}
(u,v)_S \,=\, \mathrm{Re}\,\mathrm{tr} \left(\,S^{-1}u\,S^{-1}v\,\right)
\end{equation}
which yields a Riemannian metric on $\mathcal{Q}_n$. \\
(iii) This Riemannian metric is invariant under the action of $\mathrm{GL}(n,\mathbb{H})$ on $\mathcal{Q}_n$. 
\end{proposition}
The proof of Proposition \ref{prop:metric} only requires the fact that (\ref{eq:metric1}) is a scalar product on $\mathfrak{p}_n$, and application of the rules (\ref{eq:quatrules}). It is here omitted for lack of space. 

%
%
%
%

\section{The metric in polar coordinates}
In order to provide analytic expressions in Sections \ref{sec:gauss} and \ref{sec:sampling}, we now introduce the expression of the Riemannian metric (\ref{eq:metric2}) in terms of polar coordinates. For $S \in \mathcal{Q}_n$, the polar coordinates of $S$ are the pair $(R,K)$ appearing in the decomposition (\ref{eq:spectral}). It is an abuse of language to call them coordinates, as they are not unique. However, this terminology is natural and used quite often in the literature~\cite{said1,said2}.

The expression of the metric (\ref{eq:metric2}) in terms of the polar coordinates $(R,K)$ is here given in Proposition \ref{prop:polar1}. This requires the following notation. For $i,j = 1\,,\ldots,\, n$, let $\theta_{ij}$ be the quaternion-valued differential form on $\mathrm{Sp}(n)$,
\begin{equation} \label{eq:theta1}
  \theta_{ij}(K)\,=\, \sum^n_{l=1} \,K^{\dagger}_{il}\, dK_{lj} 
\end{equation}
Note that, by differentiating the identity $K^{\dagger}K = I$, it follows that $\theta_{ij} = - \bar{\theta}_{ji}\,$. Proposition \ref{prop:polar1}  expresses the length element corresponding to the Riemannian metric (\ref{eq:metric2}).
\begin{proposition}[the metric in polar coordinates] \label{prop:polar1}
  In terms of the polar coordinates $(R,K)$, the length element corresponding to the Riemannian metric  (\ref{eq:metric2}) is given by,
\begin{equation} \label{eq:polarmetric}
   ds^2(R,K) \,=\, \sum^n_{i=1}\, dr^2_i \,+\, 8\,\, \sum_{i<j} \, \sinh^2\left(|r_i-r_j|/2\right)\,|\,\theta_{ij}|^2   
\end{equation}
where $r_i$ denote the diagonal elements of the matrix $R$. 
\end{proposition}
The proof of this proposition cannot be given here, due to lack of space.

Proposition \ref{prop:polar1} is valuable to understanding the Riemannian geometry of the space $\mathcal{Q}_n$. Precisely, it can be used to infer, with almost no calculation, the expressions of geodesics and of distance, on this space. Indeed, it becomes clear from (\ref{eq:polarmetric}) that the shortest curve connecting the identity $I \in \mathcal{Q}_n$ to a diagonal (and therefore real) element $a \in \mathcal{Q}_n$, is given by $t \mapsto a^t$ for $t \in [0,1]$. Using this simple result, and the fact that the metric (\ref{eq:metric2}) is invariant under the action of $\mathrm{GL}(n,\mathbb{H})$ on $\mathcal{Q}_n$, the equation of the minimising geodesic curve $\gamma(t)$ connecting two elements $S, Q \in \mathcal{Q}_n$ can be obtained,
\begin{equation} \label{eq:geodesic}
 \gamma(t) = \,S^{\frac{1}{2}}\,\left(\,S^{-\frac{1}{2}}QS^{-\frac{1}{2}}\,\right)^{\!\!t}\,S^{\frac{1}{2}}
\end{equation}
Accordingly, the distance between $S$ and $Q$ is
\begin{equation} \label{eq:distance}
 d(S,Q) = \left\Vert \log\left(\,S^{-\frac{1}{2}}QS^{-\frac{1}{2}}\,\right) \right\Vert
\end{equation}
where $\Vert \cdot \Vert$ is the norm corresponding to the scalar product (\ref{eq:scalproduct}).  

In (\ref{eq:geodesic}) and (\ref{eq:distance}) matrix functions, such as elevation to a power and logarithm, are computed via the decomposition (\ref{eq:spectral}), where the functions are applied to the diagonal matrix $\exp(R)$.

%
%
%
%

\section{Riemannian Gaussian distributions on $\mathcal{Q}_n$}
\label{sec:gauss}
It is possible to define Riemannian Gaussian distributions on any Riemannian symmetric space of non-positive curvature~\cite{said2}. This is indeed the case of the space $\mathcal{Q}_n\,$, as can be seen from its representation (\ref{eq:quotient}) as a quotient space, by consulting the tables  which classify irreducible Riemannian symmetric spaces of type III~\cite{besse}.

Accordingly, it is possible to define Riemannian Gaussian distributions on $\mathcal{Q}_n$. Precisely, a Riemannian Gaussian distribution on $\mathcal{Q}_n$ with Riemannian barycentre $\breve{S} \in \mathcal{Q}_n$ and dispersion parameter $\sigma > 0$ has the following probability density
\begin{equation} \label{eq:density}
p(S|\,\breve{S},\sigma) \,=\,\frac{1}{Z(\sigma)}\,\exp\left[\, -\frac{d^{\,2}(S,\,\breve{S})}{2\sigma^2}\,\right]
\end{equation}
with respect to the Riemannian volume element of $\mathcal{Q}_n$, here denoted $dv$. In this probability density, $d(S,\,\breve{S})$ is the Riemannian distance given by (\ref{eq:distance}). 

The first step to understanding this definition is computing the normalising constant $Z(\sigma)$. This is given by the integral, 
\begin{equation} \label{eq:z1}
  Z(\sigma) = \, \int_{\mathcal{Q}_n} \, \exp\left[\, -\frac{d^{\,2}(S,\,\breve{S})}{2\sigma^2}\,\right] \, dv(S)
\end{equation}
As shows in~\cite{said2}, this does not depend on $\breve{S}$, and therefore it is possible to take $\breve{S} = I$. From the decomposition (\ref{eq:spectral}) and formula (\ref{eq:distance}), it follows that 
\begin{equation} \label{eq:iddistance}
d^{\,2}(S,I) =  \sum^n_{i=1}\, r^2_i
\end{equation}
Given this simple expression, it seems reasonable to pursue the computation of the integral (\ref{eq:z1}) in polar coordinates. This is achieved in the following Proposition \ref{prop:polar2}. For the statement, write the quaternion-valued differential form $\theta_{ij}$ of (\ref{eq:theta1}) as $\theta_{ij} = \theta^a_{ij} + \theta^b_{ij}\,\mathrm{i}+\theta^c_{ij}\,\mathrm{j}+\theta^d_{ij}\,\mathrm{k}$ where $\theta^a_{ij}, \theta^b_{ij}, \theta^c_{ij}, \theta^d_{ij}$ are real-valued. 
\begin{proposition}[normalising constant] \label{prop:polar2} \\
(i)  In terms of the polar coordinates $(R,K)$, the Riemannian volume element $dv(S)$ corresponding to the Riemannian metric (\ref{eq:metric2}) is given by
\begin{equation} \label{eq:polarvolume}
   dv(R,K) =\, 8^{n(n-1)}\,\prod_{i<j} \, \sinh^4\left(|r_i-r_j|/2\right)\, \prod^n_{i=1}dr_i\, \bigwedge_{i< j} \theta^{a}_{ij}\,
  \bigwedge_{i< j} \theta^{b}_{ij}\, \bigwedge_{i< j} \theta^{c}_{ij}\, \bigwedge_{i< j} \theta^{d}_{ij}\!\!
\end{equation}
(ii) The integral $Z(\sigma)$ appearing in (\ref{eq:z1}) is given by
\begin{equation} \label{eq:z2}
  Z(\sigma) \,=\, \mathrm{Const.} \,\times\, \int_{\mathbb{R}^n}\, \exp \left(\, - \frac{1}{2\sigma^2} \,\sum^n_{i=1}\, r^2_i\,\right) \, \prod_{i<j} \, \sinh^4\left(|r_i-r_j|/2\right)\, \prod^n_{i=1}dr_i
\end{equation}
\end{proposition}
This proposition is a corollary of Proposition \ref{prop:polar1}. Formula (\ref{eq:polarvolume}) is a straightforward consequence of formula (\ref{eq:polarmetric}). Furthermore, (\ref{eq:z2}) is an immediate application of (\ref{eq:iddistance}) and (\ref{eq:polarvolume}). 

%
%
%
%

\section{Sampling and inference}
\label{sec:sampling}
The present section describes two aspects of Riemannian Gaussian distributions on $\mathcal{Q}_n$\,: i) sampling from these distributions, ii) maximum likelihood estimation of these distributions.

The first of these aspects is given in Proposition \ref{prop:sampling} below. This relies on the use of polar coordinates $(R,K)$ which appear in the decomposition (\ref{eq:spectral}). 
\begin{proposition}[Gaussian distribution in polar coordinates] \label{prop:sampling}
  Let $K$ and $r$ be independent random variables, with their values in $\mathrm{Sp}(n)$ and $\mathbb{R}^n$ respectively. Assume $K$ is uniformly distributed on $\mathrm{Sp}(n)$, and $r$ has the following probability density, with respect to the Lebesgue measure on $\mathbb{R}^n$,
\begin{equation} \label{eq:rdensity}
p(r_1\,,\ldots,\,r_n) \,\propto\,\exp \left(\, - \frac{1}{2\sigma^2} \,\sum^n_{i=1}\, r^2_i\,\right) \, \prod_{i<j} \, \sinh^4\left(|r_i-r_j|/2\right)
\end{equation}
If $S$ is given by (\ref{eq:spectral}), where the matrix $R$ has diagonal elements $r_{i\,}$, then $S$ has a Riemannian Gaussian distribution (\ref{eq:density}) with Riemannian barycentre $\breve{S} = I$ and dispersion parameter $\sigma$. Moreover, for any $\breve{S} \in \mathcal{Q}_n$ and $A \in \mathrm{GL}(n,\mathbb{H})$ such that $A\cdot I = \breve{S}$, if $Q = A\cdot S$ then $Q$ has Riemannian Gaussian distribution with Riemannian barycentre $\breve{S}$ and dispersion parameter $\sigma$. 
\end{proposition}
Proposition \ref{prop:sampling} provides a sampling algorithm for Riemannian Gaussian distributions on $\mathcal{Q}_n$. Indeed, the proposition states that in order to obtain $Q$ with Riemannian Gaussian distribution of barycentre $\breve{S}$ and dispersion $\sigma$, it is enough to know how to sample $S$ from a Riemannian Gaussian distribution with barycentre $I$. In turn, this is done using polar coordinates, through decomposition (\ref{eq:spectral}).

In this decomposition, $K$ must be sampled from a uniform distribution on $\mathrm{Sp}(n)$, and $R$ with diagonal elements $r_i$ from the multivariate density (\ref{eq:rdensity}). Sampling from a uniform distribution on $\mathrm{Sp}(n)$ can be achieved as follows\,: let $Z$ be an $n \times n$ quaternion matrix whose elements are independent normal proper quaternion random variables~\cite{nico}, and write $Z = KP$ for the polar decomposition of $Z$~\cite{zhang}. Then, $K$ has a uniform distribution on $\mathrm{Sp}(n)$. On the other hand, sampling from the multivariate density (\ref{eq:rdensity}) can be carried out using a Metropolis-Hastings algorithm, which is included in most statistical software~\cite{metropolis}. 

Consider now maximum likelihood estimation of Riemannian Gaussian distributions. This is given by the following Proposition \ref{prop:inference}. This proposition brings out the important role of the function $Z(\sigma)$ defined by (\ref{eq:z1}) and (\ref{eq:z2}). Precisely, this is the moment generating function of the Riemannian Gaussian distribution (\ref{eq:density}). If $\eta = -1/2\sigma^2$ and $\psi(\eta) = \log Z(\sigma)$, then $\psi(\eta)$ is a strictly convex function, which is the cumulant generating function of the distribution (\ref{eq:density}).
\begin{proposition}[Maximum likelihood estimation] \label{prop:inference} Let $S_1\,,\ldots,\, S_N\,$ be independent samples from a Riemannian Gaussian distribution with density (\ref{eq:density}). Based on these samples, the  maximum likelihood estimate of $\breve{S}$ is the sample Riemannian barycentre $\hat{S}_N$,
\begin{equation}
  \hat{S}_N \,=\, \mathrm{argmin}_{S \in\mathcal{Q}_n} \, \sum^N_{i=1} d^{\,2}(S_i,\,S)
\end{equation}
where the distance $d(S_i,\,S)$ is given by (\ref{eq:distance}). Moreover, the maximum likelihood estimate of $\eta = -1/2\sigma^2$ is $\hat{\eta}_N$,
\begin{equation} \label{eq:hateta}
\hat{\eta}_N = \left(\,\psi^\prime\,\right)^{-1}\left(\, \frac{1}{N} \sum^N_{i=1} d^{\,2}(S_i,\,\hat{S}_N) \,\right)
\end{equation}
where $\left(\,\psi^\prime\,\right)^{-1}$ is the reciprocal function of $\psi^\prime$, the derivative of $\psi$. 
\end{proposition}
Proposition \ref{prop:inference} indicates how the maximum likelihood estimates $\hat{S}_N$ and $\hat{\eta}_N$ can be computed. First, $\hat{S}_N$ is the sample Riemannian barycentre of $S_1\,,\ldots,\,S_N\,$. Its existence and uniqueness are guaranteed by the fact that $\mathcal{Q}_N$ is a Riemannian manifold of non-positive curvature. In practice, it can be computed using a Riemannian gradient descent algorithm~\cite{manton_globally_2004,manton_framework_2014}. Once $\hat{S}_N$ has been obtained, $\hat{\eta}_N$ is found by direct application of (\ref{eq:hateta}). This only requires knowledge of the cumulant generating function $\psi(\eta)$, which can be tabulated using the Monte Carlo method of~\cite{paolo}. 



%
%



\bibliographystyle{splncs}
\scriptsize
\bibliography{refs}

\begin{thebibliography}{10}

\bibitem{pennec}
Pennec, X.:
\newblock Intrinsic statistics on {R}iemannian manifolds: basic tools for
  geometric measurements.
\newblock J. Math. Imaging Vis. \textbf{25}(1) (2006)  127--154

\bibitem{moakher}
Chebbi, Z., Moakher, M.:
\newblock Means of {H}ermitian positive-definite matrices based on the
  log-determinant alpha-divergence function.
\newblock Linear Algebra Appl. \textbf{436}(7) (2012)  1872--1889

\bibitem{helgasson}
Helgason, S.:
\newblock Differential geometry, Lie groups, and symmetric spaces.
\newblock American Mathematical Society (2001)

\bibitem{besse}
Besse, A.L.:
\newblock Einstein manifolds, (first edition).
\newblock Springer Verlag (2007)

\bibitem{said1}
Said, S., Bombrun, L., Berthoumieu, Y., Manton, J.H.:
\newblock Riemannian {Gaussian} distributions on the space of symmetric
  positive definite matrices (accepted).
\newblock IEEE Trans. Inf. Theory (2016)

\bibitem{said2}
Said, S., Hajri, H., Bombrun, L., Vemuri, B.C.:
\newblock Gaussian distributions on {Riemannian} symmetric spaces\,:
  statistical learning with structured covariance matrices (under review).
\newblock IEEE Trans. Inf. Theory (2017)

\bibitem{flamant2016}
Flamant, J., {Le Bihan}, N., Chainais, P.:
\newblock Time-frequency analysis of bivariate signals (under review).
\newblock Applied and Computational Harmonic Analysis (2017)

\bibitem{conway}
Conway, J.H., Smith, D.A.:
\newblock On quaternions and octonions, their geometry, arithmetic and
  symmetry.
\newblock CRC Press (2003)

\bibitem{zhang}
Zhang, F.:
\newblock Quaternions and matrices of quaternions.
\newblock Linear Algebra Appl. \textbf{251} (1997)  21--57

\bibitem{kirillov}
Kirillov, A.:
\newblock An introduction to {Lie} groups and {Lie} algebras.
\newblock Cambridge University Press (2008)

\bibitem{nico}
{Le Bihan}, N.:
\newblock The geometry of proper quaternion random variables (to appear).
\newblock Signal Processing (2017)

\bibitem{metropolis}
Robert, C.P., Casella, G.:
\newblock Monte {C}arlo Statistical Methods.
\newblock Springer-Verlag (2004)

\bibitem{manton_globally_2004}
Manton, J.H.:
\newblock A globally convergent numerical algorithm for computing the centre of
  mass on compact {Lie} groups.
\newblock In: {ICARCV} 2004 8th {Control}, {Automation}, {Robotics} and
  {Vision} {Conference}, 2004. Volume~3. (December 2004)  2211--2216

\bibitem{manton_framework_2014}
Manton, J.H.:
\newblock A framework for generalising the {Newton} method and other iterative
  methods from {Euclidean} space to manifolds.
\newblock Numerische Mathematik \textbf{129}(1) (May 2014)  91--125

\bibitem{paolo}
Zanini, P., Said, S., Congedo, M., Berthoumieu, Y., Jutten, C.:
\newblock Parameter estimates of {R}iemannian {G}aussian distributions in the
  manifold of covariance matrices.
\newblock In: Sensor Array and Multichannel Signal Processsing Workshop (SAM).
  (2016)

\end{thebibliography}

\end{document}